\documentclass[11pt]{article}

\usepackage{amsmath,amssymb}
\usepackage{graphicx}
\usepackage{algorithm}
\usepackage{algpseudocode}
\usepackage{tikz}
\usepackage{longtable}
\usepackage{subcaption}
\usepackage{natbib}
\usepackage[hidelinks]{hyperref}

\usepackage[letterpaper,margin=1in]{geometry}
\usepackage[T1]{fontenc}
\usepackage{tgtermes}
\usepackage{newtxtext}
\usepackage{newtxmath}
\usepackage{bm}
\usepackage{setspace}
\bibpunct[, ]{(}{)}{,}{a}{}{,}
\def\bibfont{\small}

\usepackage{titlesec}

\titleformat{\section}
  {\normalfont\large\sffamily\bfseries}
  {\thesection.}{0.5em}{}

\titleformat{\subsection}
  {\normalfont\normalsize\sffamily\bfseries}
  {\thesubsection.}{0.5em}{}

\titleformat{\paragraph}[runin]
  {\normalfont\normalsize\bfseries}
  {}{0pt}{}[\quad]

\usepackage{caption}
\newcommand{\FIGURE}[3]{%
  #1%
  \caption{#2}%
  #3%
}

\usepackage{fancyhdr}
\begin{document}
\thispagestyle{plain} 
\begin{center}

{\LARGE\bfseries
When Solvers Get Moody: Quantifying and Predicting\\
Performance Variability in Mixed Integer Programming
\par}

\vspace{1.0em}

{\large
Mourchid Adegbindin\textsuperscript{a}
\qquad
Brenda Dietrich\textsuperscript{b}
\qquad
Jamol Pender\textsuperscript{c}
\qquad
Andrea Lodi\textsuperscript{d}
\qquad
\par}

\vspace{0.4em}

{\small
\textsuperscript{a}Cornell University, ma757@cornell.edu\\
\textsuperscript{b}Cornell University, bld34@cornell.edu\\
\textsuperscript{c}Cornell University, jjp274@cornell.edu\\
\textsuperscript{d}Cornell University, al748@cornell.edu\\
\par}

\vspace{0.8em}

{\small\itshape
Submitted to INFORMS Journal on Optimization
\par}

\end{center}

\vspace{0.5em}

\begin{abstract}
Mixed-integer programming (MIP) solvers are widely used, yet their performance can vary substantially even under seemingly neutral perturbations. Such variability complicates benchmarking, and the interpretation of computational studies. We introduce the Weighted Variability Score (WVS), a robust multi-component metric that aggregates dispersion in normalized runtime, final optimality gap, LP iterations, and branch-and-bound nodes across repeated runs. On a large subset of MIPLIB 2017 instances, each solved in 100 randomized runs under a four-hour limit, WVS reveals a right-skewed distribution with a substantial upper tail: many instances are essentially stable, but some exhibit large variability that is not explained by basic properties alone. Notably, this experiment proved optimality for the open instance \textit{neos-5045105-creuse} in two SCIP runs. Using static descriptors and root-level features, a tree-based predictor with a gated tail uplift achieves a mean absolute error of 0.058 and Spearman rank correlation of 0.65 on a sealed test split. Feature-importance analysis highlights root degeneracy, matrix structure, coefficient ranges, and presolve reductions as the main drivers of high WVS. Because variability is systematic and partially predictable from inexpensive early-run features, our findings offer a practical tool for preemptive tuning and lay a foundation for fairer evaluation of new algorithms.
\end{abstract}





\vspace{1em}



\section{Introduction}\label{sec:intro}
\begin{figure}
\FIGURE
{\centering
  \begin{subfigure}{0.49\linewidth}
    \centering
    \includegraphics[width=\linewidth]{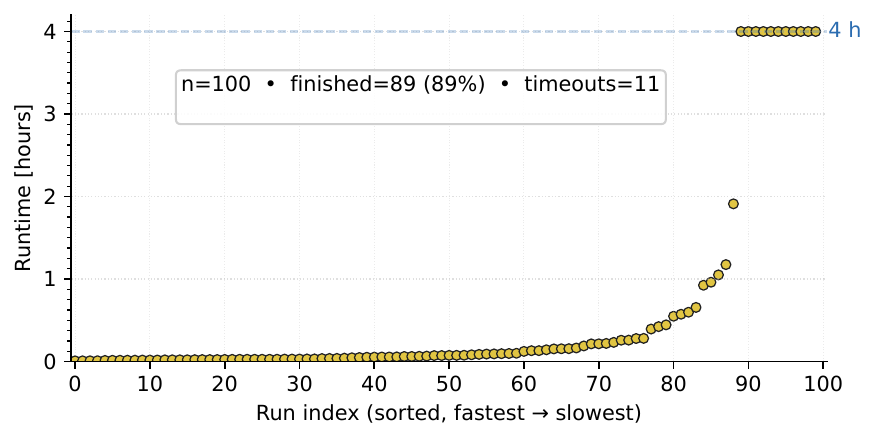}
  \end{subfigure}
  \hfill
  \begin{subfigure}{0.49\linewidth}
    \centering
    \includegraphics[width=\linewidth]{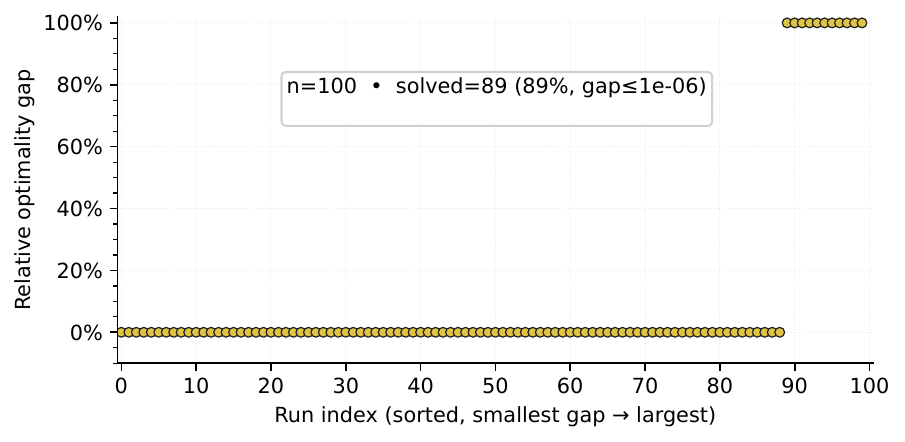}
  \end{subfigure}
}
{Empirical distributions for neos-2624317-amur. Runtime (left) and gap (right). \label{fig:fig1ab}}
{}
\end{figure}

Consider a single optimization problem, \textit{neos-2624317-amur}, from the MIPLIB~2017 benchmark library. We solve it 100 times with SCIP under identical parameter settings and a 4-hour time limit. Figure~\ref{fig:fig1ab} shows that most runs finish in well under a minute, while 11 hit the time limit at 4 hours. The total LP iterations range from $1.2\times 10^5$ to $8.1\times 10^7$, and explored branch-and-bound nodes from $2{,}837$ to $3{,}794{,}640$, even though only the random seed varies.

This behavior, where mathematically equivalent versions of the same problem, under seemingly neutral perturbations, yield different runtimes is often referred to as \textit{performance variability}. It has been discussed for decades in satisfiability and constraint programming ~\citep{GomesSelmanKautz1998_Boosting,GomesSelmanKautz2000_HeavyTails}.

\subsection{Relevant Literature}\label{subsec:intro-lit}
In mixed-integer programming (MIP), a widely cited example of performance variability is the MIPLIB instance \textit{10teams}, highlighted by \citet{danna2008performance}. Under identical settings, CPLEX solves this model on a Linux system in 0.41 seconds with 2{,}731 LP iterations and no branch-and-bound nodes, whereas on an AIX system it needs 41.39 seconds, 122{,}948 LP iterations, and 1{,}426 nodes.

Given the widespread use of MIP in industry and scientific computing \citep{Bixby2012_MIPHistory,JungerEtAl2010_50YearsIP}, performance variability has been studied in more detail. \citet{KochEtAl2011_MIPLIB2010} discuss several likely causes. Some are intrinsic to the solving process: when several branching candidates have the same scores, small perturbations in their ordering or other incidental details can lead to completely different search trees and large performance differences. Other factors include floating-point arithmetic and numerical tolerances, which can alter comparisons, change the cuts generated, and affect the bases selected by the LP solver. The same study used a permutation-based protocol: for each instance, 100 random permutations of rows and columns were generated, each permuted model was solved under identical settings, and the spread in runtimes was summarized by the coefficient of variation. Their results show that every tested case exhibits nontrivial variability, with fastest-to-slowest time ratios ranging from modest factors to several orders of magnitude.

Mittelmann's benchmark reports provide a complementary view. In earlier experiments, summarized by \citet{LodiTramontani2013}, instances in the MIPLIB~2010 \textit{benchmark set} were run once with the default seed and five times with alternative seeds, and substantial differences between default and non-default runs were observed. More recent MIPLIB~2017 reports \citep{Mittelmann2025_MIPLIB2017Bench} emphasize shifted geometric means of runtimes and solved counts under fixed limits. 
Related trajectory-based measures, such as the primal integral and the confined primal integral, summarize how solution quality evolves over time within a single run \citep{Berthold2013_PrimalIntegral,BertholdCsizmadia2021_ConfinedPI}. These are well suited for comparing progress under a fixed time budget, but they are not designed to quantify across-run dispersion on a fixed instance.
Together these studies reinforce the view that solver performance should be sampled, not inferred from a single run.

Several authors have tried either to exploit or to reduce performance variability. \textit{Bet-and-run} schemes start several short runs with different seeds and continue only the most promising one \citep{FischettiMonaci2014_BetAndRun}, and diversified parallel branch-and-cut runs play a similar role by betting on the best trajectory among different parameter settings \citep{CarvajalEtAl2014_Diversification}. Reformulation methods, including learning-based reorderings of LP and MIP models, likewise search for lucky formulations that are much easier to solve \citep{LiEtAl2023_LearningReformulate}. In parallel, there is a substantial literature on machine learning for MIP \citep{bengio2021machine} and on runtime prediction. Model-based configuration and selection methods such as ParamILS, SMAC, ISAC, and Hydra-MIP, and later runtime-prediction work, focus on choosing or tuning solvers and parameters for a given instance distribution, typically targeting mean runtime rather than variability itself \citep{HutterEtAl2009_ParamILS,HutterEtAl2011_SMAC,KadiogluEtAl2010_ISAC,XuEtAl2011_HydraMIP,HutterEtAl2014_RuntimePrediction}.

The picture that emerges is useful but incomplete. To exploit or reduce variability, one must measure it holistically and reliably for each instance and predict it from early-run information. Otherwise, sampling or stabilization schemes may be applied where they do not matter, which can be counterproductive. These observations motivate two questions: how to measure solver variability in a way that is robust and whether variability can be predicted from readily available information, in particular from features observable at or near the root node.

\subsection{Summary of contributions}\label{subsec:intro-contrib}

This paper makes four contributions. First, we introduce the weighted variability score (WVS), a robust multi-component metric. Second, we perform repeated runs on a broad MIPLIB~2017 subset with consistent limits and seed control, using SCIP on the majority of instances and Gurobi on a subset of instances that was run under the same protocol. Third, we empirically characterize the resulting variability: contributions differ across metrics, and some instances are solver-agnostic, others solver-sensitive. Finally, we build a tree-based predictor with gated tail uplift, evaluate it with tail-aware errors, and use SHAP attributions to identify key variability drivers.

The remainder of the paper is organized as follows: Section~\ref{sec:design} describes the experimental setup; Section~\ref{sec:metric} formalizes WVS; Section~\ref{sec:empirical} examines empirical patterns and solver effects; Section~\ref{sec:predict} develops and evaluates the predictive model and its interpretability; and Section~\ref{sec:discussion} discusses implications and limitations.

\section{Experimental Design}\label{sec:design}

\subsection{Benchmark Instances, solvers, and run protocol}\label{subsec:instances}

To analyze solver performance variability, we work with the sixth major release of the Mixed Integer Programming Library, MIPLIB 2017~\citep{GleixnerEtAl2021_MIPLIB2017}. The collection contains 1{,}065 instances from academic and industrial sources. All models are provided in MPS format and standardized as minimization problems. We use SCIP for its tight integration with MIPLIB and Gurobi as a commercial reference. Our goal is not to compare absolute speed but to study run-to-run variability under fixed limits.

In line with the mechanisms documented by ~\citet{danna2008performance} and ~\citet{LodiTramontani2013}, and following the practice of~\citet{KochEtAl2011_MIPLIB2010}, we perform 100 independent runs of each instance under default parameters, varying only the solver’s internal random seed (the \texttt{Seed} parameter for Gurobi and \texttt{randomseedshift} for SCIP). In SCIP, we additionally enable \texttt{permuteconss} and \texttt{permutevars} to randomly permute the constraints and variables. Each run uses a 4-hour wall-clock time limit and a common memory limit. For every run we log runtime, termination status, final relative gap, branch-and-bound node count, total LP iterations, and a collection of root-node statistics described in Section~\ref{subsec:features}. Instances with more than 10 incomplete runs are excluded. After this filtering, 1009 instances remain and form our main corpus. Of these, 388 instances were also run with Gurobi under the same protocol; this common set is the subset used in the cross-solver analyses.

All experiments are executed on identical Microsoft Azure Standard FX48mds virtual machines with Intel Xeon Gold processors and a fixed amount of RAM per core~\citep{AzureFX48mds}. Autoscaling and background jobs are disabled, and CPU frequency scaling is constrained as much as the platform allows. All solves are forced to single-threaded mode. Experiments use SCIP~9.2.1 with SoPlex~7.1.3 (accessed via PySCIPOpt~5.4.1)~\citep{BestuzhevaEtAl2023_SCIPSuite,MaherEtAl2016_PySCIPOpt} and, on the 388-instance common subset described above, Gurobi~12.0.2~\citep{Gurobi2025_Manual}. Although cloud infrastructure cannot eliminate all exogenous timing noise, the execution environment was kept as uniform as the platform allows. Since WVS is not based on runtime alone, incidental wall-clock fluctuations do not determine the conclusions.

\subsection{Feature Extraction}\label{subsec:features}
Our features can be divided into two groups: static structural features, which are fixed for the instance, and early-run dynamic features, which summarize inexpensive information from presolve and the root node. Static features describe the scale and basic structure of the model. They include numbers of variables, constraints, and nonzeros; variable-type ratios (binary, general-integer, continuous); simple measures of sparsity and row/column aspect ratio; and summary statistics of objective coefficients and right-hand sides. We also include a small number of structural flags such as counts of set-partitioning, set-covering, knapsack, indicator, precedence, and cardinality constraints; counts of bounded variables, following the MIPLIB~2017 tag collection~\citep{MIPLIB2017_TagCollection}. Dynamic features are extracted from SCIP runs and capture early presolving and root-node behavior under the standardized run protocol. For each instance we record basic timing and effort at presolve and at the root LP: total presolving time, the number of presolving rounds at different effort levels (fast, medium, exhaustive), the numbers of variables, constraints, and nonzeros remaining after presolve, root LP solving time, and the number of LP iterations. We then derive several degeneracy and dual-activity proxies. Additional root-level signals include the number of fractional variables, the number and type of branching and pseudo-branching candidates, and symmetry indicators.

Some of these quantities, such as the average number of simplex pivots before a significant improvement, are, to our knowledge, not standard in the empirical MIP literature and were designed specifically for this study to capture aspects of root behavior that may relate to downstream variability. When a dynamic quantity can vary across seeds, we summarize it at the instance level by a simple aggregation (for example, a mean or median over the 100 runs).

The complete run-level data and code used to reproduce the computational results will be made publicly available upon publication.
\section{Defining and Measuring Solver Variability}\label{sec:metric}

\subsection{Motivation}
\label{subsec:metric-motivation}

Most quantitative measures of performance variability in MIP are time-based. In the MIPLIB~2010 study, \citet{KochEtAl2011_MIPLIB2010} summarize time-to-optimality dispersion via a variability score (VS) for each instance, defined as the sample standard deviation of solution times divided by their mean (a coefficient of variation). This statistic is scale-free and easy to compute, but it ignores solution quality and work measures, and it is defined on time-to-optimality, so timeouts and unsolved runs require ad hoc handling. Cross-instance benchmarking tools such as performance profiles \citep{DolanMore2002} and shifted geometric means, as used in Mittelmann’s reports, are designed to compare solvers across collections of instances, not to provide a per-instance variability score.

For our purposes we need a \emph{per-instance} index of solver instability that (i) reflects dispersion in several performance indicators (time, gap, work), (ii) is scale-free, (iii) handles timeouts naturally, (iv) is not dominated by a few extreme runs, and (v) remains simple enough to explain and to decompose back into component contributions. The construction below starts from four per-run components and builds a modified coefficient of variation (MCV) for each, then combines them into a single scalar index $\mathrm{WVS}$. We refer to runtime, final gap, LP iterations, and node counts as performance measures. When these are aggregated inside $\mathrm{WVS}$, we sometimes refer to their respective contributions as the four components of the score. 

\subsection{Components and per-component dispersion}
\label{subsec:metric-components}

For each instance $i$, solver $s$, and run $j = 1,\dots,N$ we record four basic performance measures.

\begin{enumerate}
\item \textbf{Normalized time $T_{i,s,j}$.}
Let $T$ be the wall-clock time limit, and $\tilde t_j$ the observed runtime of run $j$, capped at $T$. We define
\[
t_j \;=\; \frac{\max\{1,\tilde t_j\}}{T} \in (0,1],
\]
so that all times are normalized to $(0,1]$ and very small runtimes are rounded up to one second to avoid numerical artifacts. Time-limited runs have $t_j = 1$. We write $T_{i,s,j} = t_j$.

\item \textbf{Final relative gap $G_{i,s,j}$.}
Let $\tilde g_j$ denote the final relative gap reported by the solver. We cap this quantity at one,
\[
g_j \;=\; \min\{1,\tilde g_j\},
\]
and set $G_{i,s,j} = g_j$.

\item \textbf{LP iterations $I_{i,s,j}$.}
Total LP iterations over the entire run.

\item \textbf{Branch-and-bound nodes $B_{i,s,j}$.}
Total processed branch-and-bound nodes.
\end{enumerate}

For each component $c \in \{T,G,I,B\}$ and instance–solver pair $(i,s)$ this gives a sample $
x^{(c)}_1,\dots,x^{(c)}_N$
of length $N$, one value per seed. A natural scale-free dispersion measure is the classical coefficient of variation $s/\mu$, where $\mu$ is the empirical mean and
\[
s \;=\; \sqrt{\frac{1}{N}\sum_{j=1}^{N}(x_j-\mu)^2}
\]
is the empirical standard deviation. In our setting, however, this measure reacts strongly to a few extreme runs and can become arbitrarily large when $\mu$ is small; see, for example, \citealp{Allison1978_Inequality,ArachchigeEtAl2022_RobustCV}.

For each component we define a modified coefficient of variation (MCV). 
For a fixed instance--solver pair, let $x_1,\dots,x_N$ be the $N$ run 
values with empirical mean $\mu$ and standard deviation $s$ as defined above. The per-component 
MCV is
\begin{equation}
\mathrm{MCV}_c =
\begin{cases}
\displaystyle
\sqrt{\dfrac{\sum_{j=1}^N w_j (x_j - \mu)^2}{\mu^2 \sum_{j=1}^N w_j}},
& \mu > 0,\\[6pt]
0, & \mu = 0,
\end{cases}
\label{eq:mcv}
\end{equation}
where $w_j = 1/(1+z_j^2)$ and $z_j = (x_j - \mu)/s$ if $s>0$ (otherwise
$z_j = 0$).

This is a scale-free analogue of the classical coefficient of variation that down-weights extreme standardized deviations, so a few outlying runs have limited influence.

\subsection{The Weighted Variability Score (WVS)}
\label{subsec:metric-wvs}

For a given instance--solver pair $(i,s)$, let
$\mathrm{MCV}^T_{i,s},\quad \mathrm{MCV}^G_{i,s},\quad
\mathrm{MCV}^I_{i,s},\quad \mathrm{MCV}^B_{i,s}
$
denote the modified coefficients of variation for normalized time, gap, iterations, and branch-and-bound nodes. For any weight vector $
w = (w_T, w_G, w_I, w_B) $ 
with nonnegative entries summing to one, we define the \emph{Weighted Variability Score} (WVS) as
\begin{equation}
\mathrm{WVS}_{i,s}(w)
\;=\;
w_T\,\mathrm{MCV}^T_{i,s}
\;+\;
w_G\,\mathrm{MCV}^G_{i,s}
\;+\;
w_I\,\mathrm{MCV}^I_{i,s}
\;+\;
w_B\,\mathrm{MCV}^B_{i,s}.
\label{eq:WVS_Linear_Comb}
\end{equation}

Time and final gap reflect what users observe most directly, while iterations and nodes act as proxies for search effort. We choose $w$ on the SCIP--Gurobi overlap set by searching over a
discrete grid on the simplex
\[
\Delta=\{w\in\mathbb{R}_{\geq 0}^{4}:
w_T+w_G+w_I+w_B=1\}.
\]
For each candidate weight vector, we count the number of overlap
instances for which the absolute cross-solver disagreement exceeds a
fixed calibration threshold $\tau_w = 0.20$:
\[
\min_{w\in\Delta}
\sum_{i\in\mathcal I_{\mathrm{overlap}}}
\mathbf{1}
\left\{
\left|
\mathrm{WVS}_{i,\mathrm{SCIP}}(w)
-
\mathrm{WVS}_{i,\mathrm{Gurobi}}(w)
\right|>\tau_w
\right\}.
\]
Ties are broken by lower mean absolute disagreement and then by higher Spearman rank correlation. Several nearby weight vectors perform similarly; as a representative choice we fix $(w_T,w_G,w_I,w_B)=(0.65,0.05,0.20,0.10)$.
Alternative solutions of the calibration problem and sensitivity plots induce very similar per-instance rankings. For instances solved by both solvers we also form a \emph{consensus} label. Let $
\Delta_i \;=\; \bigl|\mathrm{WVS}_{i,\mathrm{SCIP}} - \mathrm{WVS}_{i,\mathrm{Gurobi}}\bigr|$
under the fixed weights above. Given a tolerance $\tau_{\mathrm{gs}}$, we define
\begin{equation}
\mathrm{WVS}^\ast_i \;=\; \operatorname{median}\{\mathrm{WVS}_{i,\mathrm{SCIP}}, \mathrm{WVS}_{i,\mathrm{Gurobi}}\}
\quad\text{whenever}\quad \Delta_i \le \tau_{\mathrm{gs}}.
\label{eq:WVS_gurobi_Scip}
\end{equation}
Overlap instances with $\Delta_i > \tau_{\mathrm{gs}}$ are marked as solver-discordant and excluded from supervised learning but retained for descriptive analysis. For all other instances (SCIP-only), we simply take $\mathrm{WVS}^\ast_i = \mathrm{WVS}_{i,\mathrm{SCIP}}$. In the experiments below, we set \(\tau_{\mathrm{gs}}=0.15\).

By construction, WVS is multi-component, scale-free, and timeout-robust, and it can be decomposed back into its four components through the $\mathrm{MCV}$ values. We find it useful to think of WVS as a variability ``condition number'' for an instance--solver pair: it summarizes typical spread in time, quality, and work, is stable under a few extreme runs, and remains simple enough to interpret in terms of its underlying components.

\section{Empirical Characterization of Variability}\label{sec:empirical}
\subsection{Overall Variability}\label{subsec:emp-overall}
\begin{figure}
\FIGURE
{\centering
  \begin{subfigure}[b]{0.37\linewidth}
    \centering
    \includegraphics[width=1.1\linewidth]{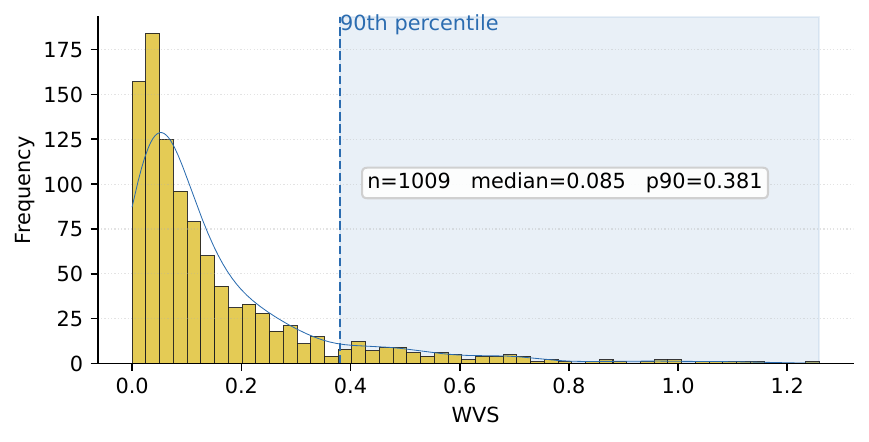}
    \caption{Distribution of $WVS_{i,\text{SCIP}}$}
    \label{fig:wvs_overview-dist}
  \end{subfigure}
  \hfill
  \begin{subfigure}[b]{0.3\linewidth}
    \centering
    \includegraphics[width=.9\linewidth]{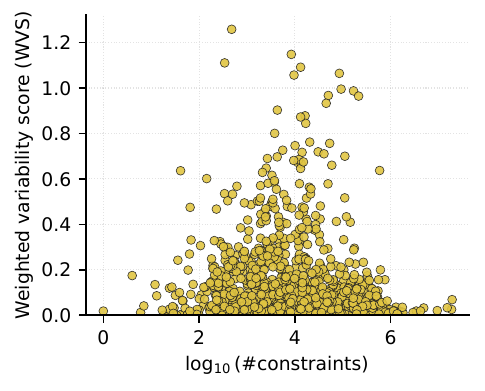}
    \caption{$\mathrm{WVS}_{i,\text{SCIP}}$ vs.\ constraints}
    \label{fig:wvs_overview-size}
  \end{subfigure}
  \hfill
  \begin{subfigure}[b]{0.3\linewidth}
    \centering
    \includegraphics[width=.9\linewidth]{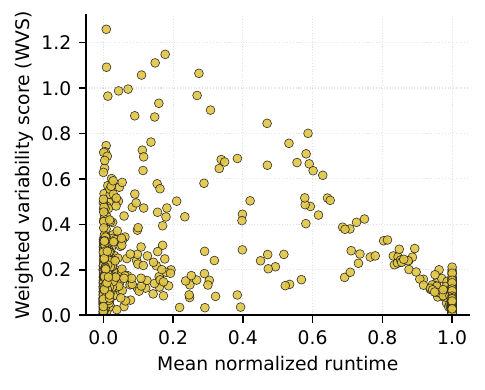}
    \caption{$\mathrm{WVS}_{i,\text{SCIP}}$ vs.\ runtime}
    \label{fig:wvs_overview-time}
  \end{subfigure}
}
{Overall behavior of the variability score $\mathrm{WVS}_{i,\text{SCIP}}$. In panels~(b)--(c) each point is one instance; the wide vertical spread
shows that neither size nor runtime alone explains variability.\label{fig:wvs_overview}}
{}
\end{figure}

Unless noted otherwise, statistics in this section use the per-instance variability
score $\mathrm{WVS}_{i,\text{SCIP}}$ based on $100$ randomized runs with SCIP.
When we distinguish solvers, we write $\mathrm{WVS}_{i,\text{SCIP}}$ and
$\mathrm{WVS}_{i,\text{Gurobi}}$.  Figure~\ref{fig:wvs_overview}\subref{fig:wvs_overview-dist} shows the empirical distribution of $\mathrm{WVS}_i$ across the corpus. The distribution is clearly right-skewed: many instances exhibit modest variability, while a non-negligible minority occupies a higher-variability regime separated from the bulk by a visible shoulder. We use the empirical 90th percentile of $\mathrm{WVS}_i$ as a simple ``high-variability'' threshold in later analyses.

A natural question is whether $\mathrm{WVS}_{i,\text{SCIP}}$ is simply
repackaging instance size or baseline difficulty.  Panels~\ref{fig:wvs_overview-size}
and~\ref{fig:wvs_overview-time} plot $\mathrm{WVS}_{i,\text{SCIP}}$ against
$\log_{10}$ of the number of constraints and against the \emph{mean}
normalized runtime across seeds. Variability does increase somewhat among
very large or long-running instances, as one would expect, but high
$\mathrm{WVS}_i$ values appear across the full size and runtime ranges:
some relatively small, fast models show pronounced variability, while
some large or slow models have very small $\mathrm{WVS}_i$ because all
runs behave similarly. Correlations between $\mathrm{WVS}_i$ and these
simple proxies are moderate rather than near one, so $WVS_i$ behaves as
an additional axis, complementary to size and baseline difficulty.

\paragraph{An Open Instance Solved}\label{subsec:open-instance-solved}

One instance provides a particularly striking example of the practical consequences of variability. At the time of our experiments, \textit{neos-5045105-creuse} was listed as open in MIPLIB~2017. Out of the 100 SCIP runs, exactly two proved optimality within the four-hour limit. The two successful runs were themselves very different. One terminated in just 122.77 seconds, 101{,}502 LP iterations, and 15{,}517 branch-and-bound nodes, whereas the other required 10{,}424.45 seconds, 17{,}716{,}828 LP iterations, and 3{,}172{,}681 nodes. This case demonstrates that performance variability can occasionally yield serendipitous breakthroughs on difficult instances.

\subsection{Component-Level Behavior}\label{subsec:emp-components}

\begin{figure}[htbp]
\centering

\vspace{1em} 

\begin{subfigure}{\linewidth}
  \centering
  \includegraphics[width=.6\linewidth]{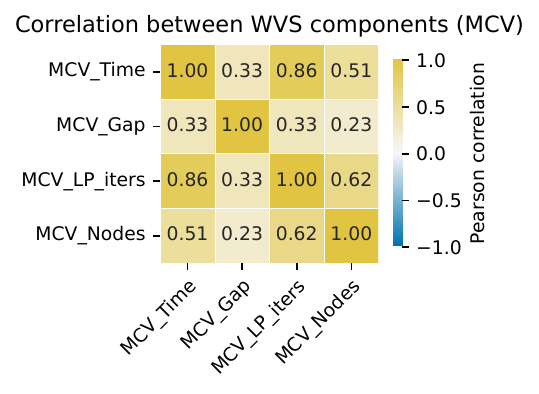}
  \caption{Correlations among component dispersions.}
  \label{fig:mcv_corr}
\end{subfigure}\\[0.25em]

\begin{subfigure}{\linewidth}
  \centering
  \includegraphics[width=.6\linewidth]{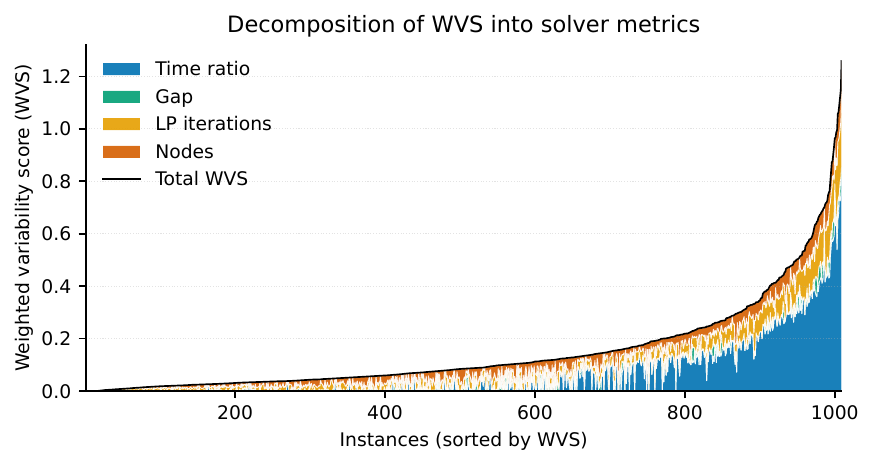}
  \caption{Component contributions across instances sorted by WVS.}
  \label{fig:mcv_deciles}
\end{subfigure}\\[0.25em]

\begin{subfigure}{\linewidth}
  \centering
  \includegraphics[width=.6\linewidth]{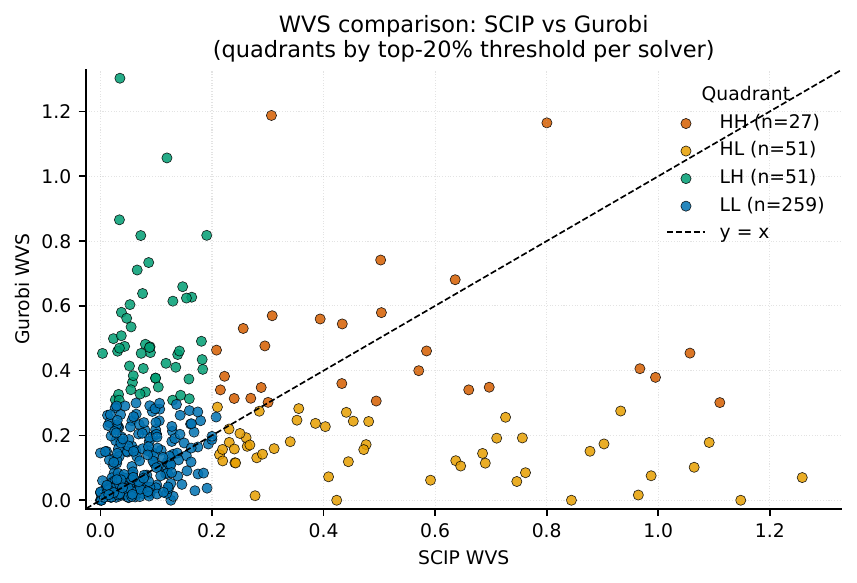}
  \caption{Cross-solver comparison of variability.}
  \label{fig:cross_solver}
\end{subfigure}
\caption{Component-level behavior and cross-solver comparison.}
\label{fig:components}
\end{figure}

The score $\mathrm{WVS}$ aggregates four components: dispersion in
normalized time, final gap, LP iterations, and node counts.
Figure~\ref{fig:components}\subref{fig:mcv_corr} shows the correlation
matrix of the component-wise dispersions $\mathrm{MCV}^T$,
$\mathrm{MCV}^G$, $\mathrm{MCV}^I$, and $\mathrm{MCV}^B$ across
instances for SCIP. LP-iteration dispersion tracks time dispersion
closely (correlation about $0.86$), as expected, but remains useful on
instances where many runs hit the time limit and time is censored.
Correlations involving gap and nodes are more modest, so no single
component fully explains the others.

To see how the four components contribute along the variability
spectrum, we sort instances by $\mathrm{WVS}_i$ and decompose
$\mathrm{WVS}_{i,\text{SCIP}}$ into its weighted contributions from
time, gap, iterations, and nodes. Panel~\ref{fig:mcv_deciles} plots
these contributions. Near the low-variability end, all contributions
are small. Toward the high-variability end, the time contribution
grows sharply and becomes dominant on average, while iterations and
nodes continue to contribute in an instance-dependent way. The gap
contribution is generally smaller under the final weighting.

\subsection{Historical Consistency with MIPLIB 2010}\label{subsec:emp-historical}
\begin{figure}
\FIGURE
{\centering
\includegraphics[width=.65\linewidth]{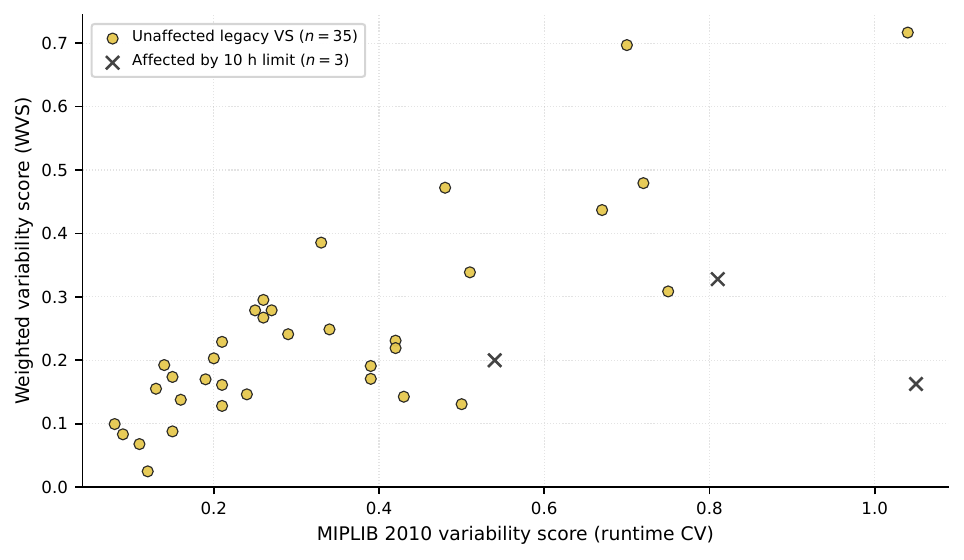}
}
{MIPLIB~2010 runtime-based variability score (VS) versus current SCIP WVS on 38 common instances. Crosses mark three legacy VS values reported as affected by the 10-hour limit. Spearman $\rho=0.675$ overall and $0.740$ when these three observations are excluded.
\label{fig:miplib2010-vs-wvs}}
{}
\end{figure}

As a historical consistency check, we compare WVS with the runtime-based variability score (VS) reported in the MIPLIB~2010 study~\citep{KochEtAl2011_MIPLIB2010}. Of the instances with a published numeric VS, 38 also appear in our 1{,}009-instance SCIP corpus. Figure~\ref{fig:miplib2010-vs-wvs} plots the two scores on this common subset. Despite differences in solver generation, hardware, perturbation protocol, and score construction, the two rankings show a substantial positive association (Spearman $\rho=0.675$). As expected for comparison with a runtime-only historical metric, the ordinary runtime coefficient of variation from our current runs is slightly more correlated with the legacy VS ($\rho=0.680$). Thus, WVS retains much of the historical runtime-variability signal while also incorporating dispersion in final gap and search effort. Three legacy VS values were reported as affected by the 10-hour limit in the MIPLIB~2010 experiment; excluding them as a sensitivity check gives $\rho=0.740$ for WVS and $\rho=0.767$ for the current runtime CV.

\subsection{Cross-Solver Alignment and Discordant Instances}\label{subsec:emp-cross}

We now examine how variability behaves across solvers on the overlap set
of $388$ instances where both SCIP and Gurobi were run.
Figure~\ref{fig:cross_solver} plots $\mathrm{WVS}_{i,\text{SCIP}}$
against $\mathrm{WVS}_{i,\text{Gurobi}}$ for these instances. Most mass
lies near the diagonal: The scatter shows a moderate positive association across solvers, together with substantial solver-specific variation. Rank-based correlation confirms this overall alignment (Spearman $\rho=0.403$). To highlight solver-sensitive cases, we define per-solver high-variability thresholds as the 80th percentiles of $\mathrm{WVS}_{\cdot,\text{SCIP}}$ and $\mathrm{WVS}_{\cdot,\text{Gurobi}}$ on the overlap set and classify each instance into four quadrants: HH (high for both), LL (low for both), HL (high for SCIP only), and LH (high for Gurobi only). In Figure~\ref{fig:cross_solver}, points are colored by these labels. Most instances fall into LL; a smaller but nontrivial fraction lies in HH, suggesting that their instability is largely instance-driven. The remaining HL and LH cases form a discordant subset where one solver is stable and the other highly variable, which is particularly relevant for algorithm design and motivates the consensus-label construction in Section~\ref{sec:metric}.

\section{Predictive Modeling of Variability}\label{sec:predict}
\subsection{Problem Setup and Evaluation Metrics}\label{subsec:pred-setup}

For each instance $i$ we start from the static and dynamic features defined in Section~\ref{subsec:features}, extracted once from the solver logs, and collect them into a feature vector $x_i$. Static features describe scale and structure (sizes, type ratios, sparsity, coefficient statistics, presolve summaries); dynamic features, available for SCIP, capture early presolve and root-node behavior under the standardized protocol (timings, iteration counts, degeneracy and dual-activity proxies, branching candidates, symmetry and decomposition indicators, ``pivots before improvement,'' and related quantities).

\paragraph{Feature engineering and screening.}
Before modeling we apply light but systematic feature engineering: skewed counts and effort measures are mapped with $\log(1+x)$, several base quantities are converted to ratios or clipped to safe ranges, and we add a small set of interaction and indicator features aimed at amplifying tail behavior (e.g., products of degeneracy and gap proxies, degeneracy interacted with tightened bounds, row/column imbalance flags, and binary- or integer-heavy indicators). On the training split only, we remove constant or near-constant predictors and greedily prune correlated ones using a shallow tree–based importance ranking, yielding a single global shortlist of roughly one hundred features that is fixed before any test-time evaluation and shared by all models.

\paragraph{Prediction task and evaluation.}
The prediction task is to learn a mapping from $x_i$ to a predicted value $\hat y_i$ that approximates the consensus variability label $WVS_i^\ast$ defined in Section~\ref{sec:metric}, using these structural and instance-level aggregated early-run features. After applying the consensus-label and modeling filters, 820 instances remain. We create a single frozen $90{:}10$ train--test split stratified by deciles of $WVS^\ast$, so that the test set preserves the prevalence of high-variability instances. All model selection and tuning use $k$-fold cross-validation inside the training portion; the test set is untouched until the final evaluation. Where indicated, we allow mild oversampling of mid- and high-$WVS^\ast$ strata within the training split, under a global cap; validation and test folds are never replicated. We evaluate prediction quality with both global and tail-aware metrics. On the test set we report
(i) global mean absolute error (MAE) between $\hat y_i$ and $WVS_i^\ast$ and Spearman rank correlation $\rho$ over all instances, and
(ii) MAE and Spearman $\rho$ restricted to the subset of instances in the top decile of \emph{true} $WVS_i^\ast$.

\subsection{Models: Baselines and Final Predictor}\label{subsec:pred-models}

We compare three models on the same screened feature set: (i) linear baselines (ordinary least squares and ridge on z-scored features); (ii) a gradient-boosted decision tree regressor \citep{Friedman2001_GradientBoosting} with a Huber loss \citep{Huber1964_RobustLocation}, tuned by cross-validation; and (iii) a slightly richer gated-expert variant of this tree model. We also tried other tabular learners (Elastic Net, alternative boosted trees, a shallow neural network), but none improved tail MAE or Spearman's $\rho$, so we use the Huber boosted tree as the main baseline.

\paragraph{Final predictor: gated experts with nonnegative uplift}

The final predictor uses the same screened feature set and CatBoost
framework but adds simple specialization for the high-variability
regime. We fit three tree ``experts'' on the training data
\citep{ProkhorenkovaEtAl2018_CatBoost}: a \emph{base} expert with a
Huber loss over the full range of $WVS^\ast$, a \emph{synthetic-tail}
expert trained with the same loss but with SMOTE-style synthetic
oversampling of mid- and high-$WVS^\ast$ training instances
\citep{ChawlaEtAl2002_SMOTE}, and a \emph{tail} expert trained with a
$0.90$ quantile loss \citep{KoenkerBassett1978_RegressionQuantiles}
focused on the upper range. Using cross-validation, we obtain
out-of-fold predictions $\hat y_i^{\text{base}}$,
$\hat y_i^{\text{sm}}$, and $\hat y_i^{\text{q90}}$ on the training
instances.

Treating $\hat y_i^{\text{base}}$ as a one-dimensional index, we define
two logistic gates $g_i^{\text{sm}}, g_i^{\text{q90}} \in [0,1]$,
with thresholds determined from out-of-fold training predictions, so
that the specialized experts are gradually activated as
$\hat y_i^{\text{base}}$ moves into the upper range. To keep the model
conservative, the specialized experts are only allowed to add
nonnegative corrections. We define the capped uplifts
\[
u_i^{\text{sm}}
=
\min\!\left\{0.08,\,
\max\!\left(\hat y_i^{\text{sm}}-\hat y_i^{\text{base}},0\right)
\right\},
\]
and
\[
u_i^{\text{q90}}
=
\min\!\left\{0.08,\,
\max\!\left(\hat y_i^{\text{q90}}-\hat y_i^{\text{base}},0\right)
\right\}.
\]
We then fit nonnegative coefficients
$\beta_{\text{sm}},\beta_{\text{q90}}$ by nonnegative least squares
using the out-of-fold training predictions. The final predictor has
the form
\begin{equation}
\hat y_i
=
\hat y_i^{\text{base}}
+
\beta_{\text{sm}}\, g_i^{\text{sm}}\, u_i^{\text{sm}}
+
\beta_{\text{q90}}\, g_i^{\text{q90}}\, u_i^{\text{q90}},
\label{eq:final-predictor}
\end{equation}
with gates and coefficients frozen before any test evaluation.

\subsection{Training Protocol and Main Results}\label{subsec:pred-results}

\begin{figure}
\FIGURE
{\centering
  \includegraphics[width=.9\linewidth]{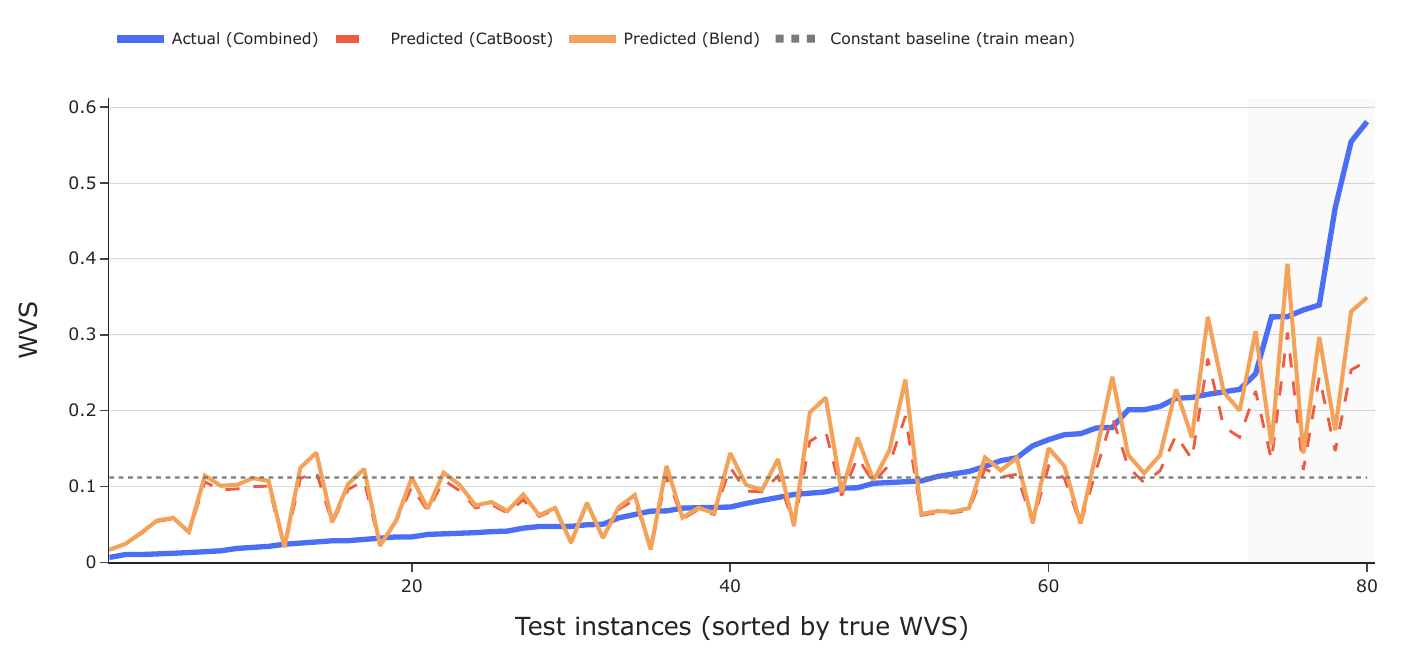}
}
{Sorted true vs.\ predicted $WVS^\ast_i$. The blended model achieves a mean absolute error (MAE) of $0.058$ and a Spearman rank correlation of $0.65$.\label{fig:pred-sorted}}
{}
\end{figure}

Tree hyperparameters are tuned by cross-validation on the training split,
balancing global MAE and MAE in the top decile of $WVS^\ast$; the uplift
coefficients are fit from the out-of-fold training predictions as described above. Figure~\ref{fig:pred-sorted} shows true and predicted $WVS_i^\ast$ on the sealed test set, with instances sorted by the true value and the top decile shaded.

We use the training protocol described in Section~\ref{subsec:pred-setup} to tune tree depth, learning rate, number of boosting iterations, and regularization parameters, balancing global MAE and MAE on the top decile of $WVS^\ast$. On the sealed test split, the final gated uplift model achieves a mean absolute error of about $0.058$ and a Spearman rank correlation of about $0.65$ with $WVS^\ast$. The Huber-tree baseline already improves substantially over ridge regression, and adding the tail-focused expert further sharpens rank quality and reduces error on the most variable instances without degrading global fit. As a robustness check, we also trained the same pipeline separately for each component dispersion (normalized time, final gap, LP iterations, nodes) and then formed a composite predictor by a nonnegative linear combination fit, on training predictions, to the fixed weights used in the definition of $WVS$. The resulting ``sum-of-components'' model attains test MAE and rank correlation close to those of the direct $WVS^\ast$ head.

\begin{figure}[htbp]
\centering
\caption{SHAP-based interpretability for the final predictor.}
\label{fig:shap}
\vspace{1em} 

\begin{subfigure}{\linewidth}
  \centering
  \includegraphics[width=.7\linewidth]{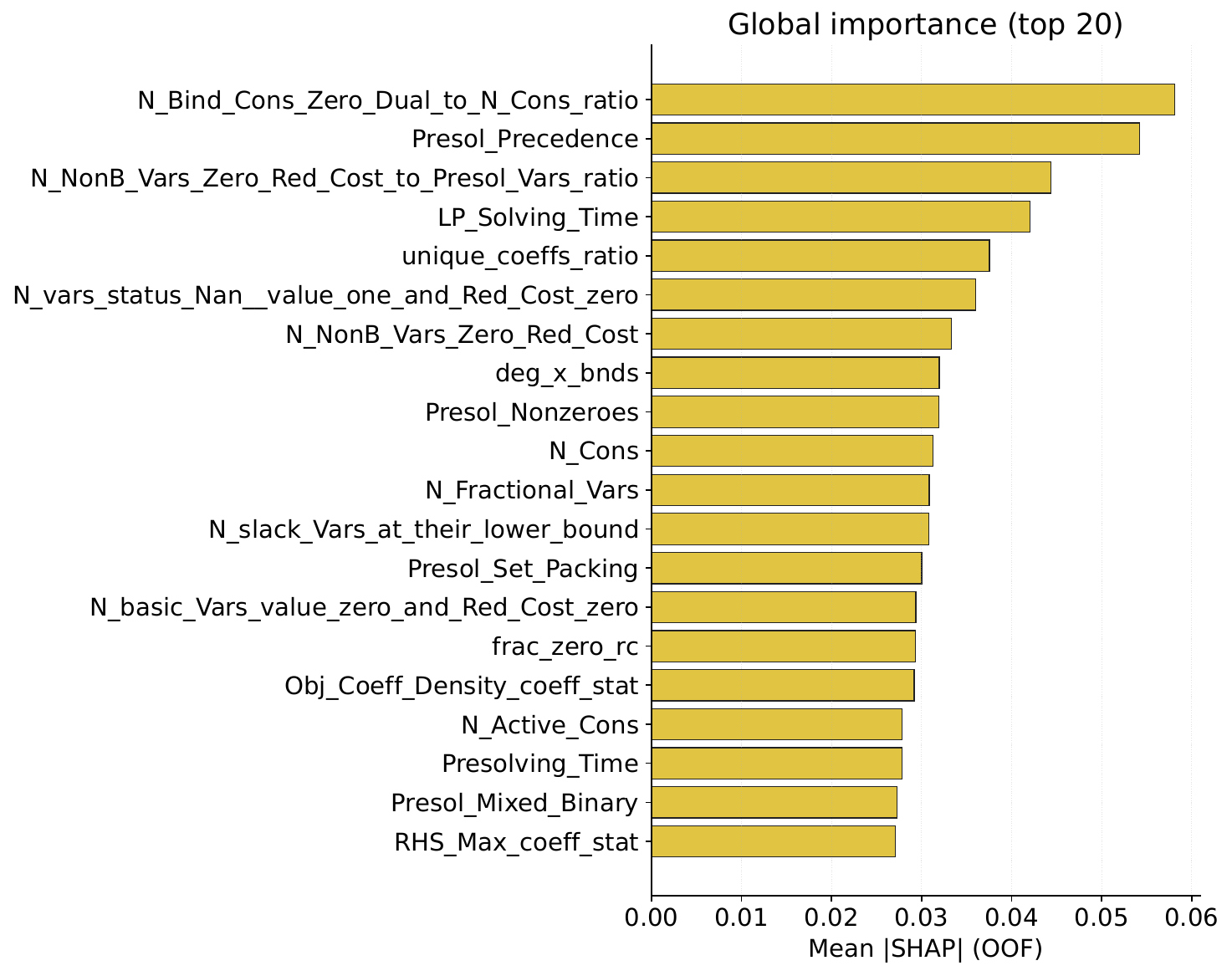}
  \caption{Top $20$ features by mean abs SHAP value.}
  \label{fig:shap-global}
\end{subfigure}\\[0.2em]

\begin{subfigure}{\linewidth}
  \centering
  \includegraphics[width=.9\linewidth]{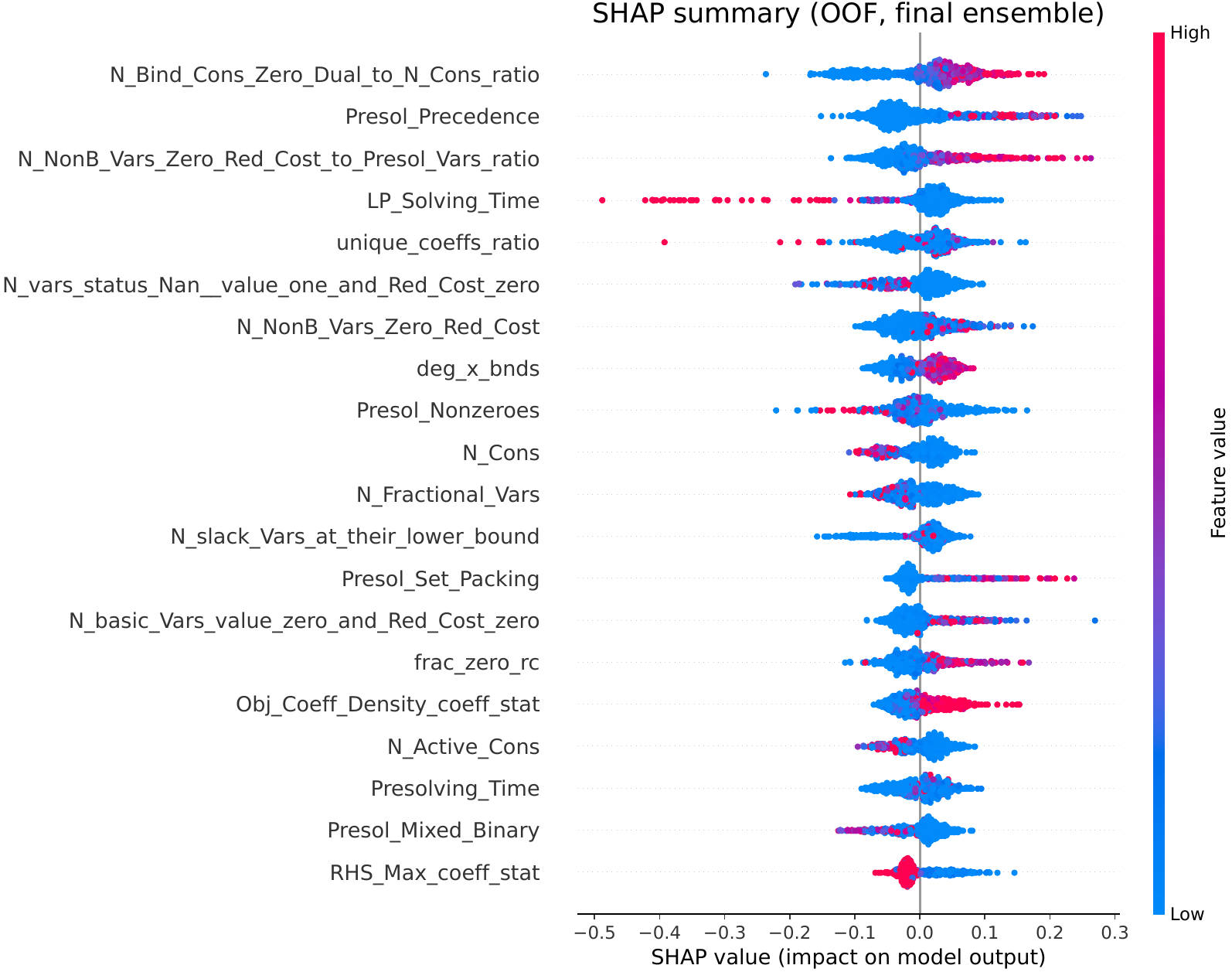}
  \caption{Beeswarm plot for the same features.}
  \label{fig:shap-beeswarm}
\end{subfigure}
\end{figure}

\paragraph{Interpretability}\label{subsec:pred-interpret}
To interpret the logic behind the model’s predictions, we compute TreeSHAP values for the base and tail experts. By combining these values according to the gating and uplift logic in~\eqref{eq:final-predictor}, we obtain a per-instance attribution vector that decomposes the final prediction into individual feature contributions. Figure~\ref{fig:shap} summarizes these values across instances. Panel~\ref{fig:shap-global} ranks the top 20 features by mean absolute SHAP value, which measures the average magnitude of a feature’s contribution to the prediction. Panel~\ref{fig:shap-beeswarm} provides a more granular view through a beeswarm plot, which shows both the direction and the spread of these effects. In the beeswarm plot, each point represents a single MIP instance. Its position on the horizontal axis is the SHAP value: points to the right of zero increase the predicted $WVS^\ast$, while points to the left decrease it. The color of each point represents the feature value, from low (blue) to high (red). For example, the ratio of binding constraints with zero dual values to the total number of constraints shows a clear monotone pattern: high values tend to push the prediction upward, consistent with the role of dual degeneracy in solver instability. Where colors overlap near zero, the effect is weaker or more contingent on interactions with other descriptors. The broader implications of these patterns are discussed in Section~\ref{sec:discussion}.

\section{Discussion}\label{sec:discussion}

\paragraph{What Drives Variability}\label{subsec:disc-drivers}

It is incomplete to describe performance variability by time alone. Two runs of the same instance can have similar runtimes but very different behavior in solution quality and solver effort, and when many runs time out we still observe variability in gaps, iterations, and nodes. The aim of $WVS$ is therefore to treat time, quality, and work jointly. The empirical results confirm this: runtime-only proxies correlate with $WVS^\ast$ but do not fully explain it, especially in the upper tail, and many high-$WVS^\ast$ cases are driven as much by gap and effort dispersion as by time.

The component-wise models sharpen this picture. Training the same pipeline separately on dispersion in normalized time, final gap, LP iterations, and nodes shows that time is the easiest component to predict, while gap and work components are noisier but clearly not random. A composite built from the four component heads attains test performance very close to that of the direct $WVS^\ast$ head, indicating that the index is not dominated by a single dimension: volatile instances are those where time, gap, and effort all vary meaningfully across runs. Consistently, runtime-only variability proxies underperform the $WVS^\ast$ predictor when the goal is to identify the top-$WVS^\ast$ instances.

The SHAP analysis of the final predictor \citep{LundbergLee2017_SHAP,LundbergEtAl2020_TreeExplainer}, implemented via scikit-learn-style pipelines \citep{PedregosaEtAl2011_ScikitLearn}, explains which mechanisms drive this behavior. It shows that $WVS^\ast$ is driven by a combination of root-level degeneracy, coefficient patterns, and the combinatorial structure exposed during presolve, rather than by instance scale alone. The most consistent signals of instability are indicators of primal and dual inactivity: the ratio of binding constraints with zero dual values and the fraction of non-basic variables with zero reduced costs. These are classic signatures of a degenerate root relaxation in which many nearly equivalent bases create flat directions in the search space. In such settings, small perturbations in pricing or pivoting can propagate into substantially different search trees, which explains why two randomized runs of the same instance may exhibit high dispersion.

The SHAP analysis also highlights that variability depends on the specific structure remaining after presolve. Features such as the counts of precedence constraints, mixed-binary variables, and active constraints after presolve all rank highly. Notably, the beeswarm plot (Figure~\ref{fig:shap-beeswarm}) suggests that large presolved cores (high `Presol\_Nonzeros`) and long root LP-solving times are not, by themselves, reliable indicators of instability once other structural features are taken into account. Variability is therefore not a simple monotone function of problem size or computational effort.

The role of numerical uniqueness further supports this view. A high `unique\_coeffs\_ratio`, the ratio of unique objective coefficients to total non-zeros, tends to be associated with lower predicted variability, whereas repetitive coefficient patterns are associated with higher predicted variability. This is consistent with the hypothesis that repeated coefficients create more cost-based ties, increasing the opportunities for solver decisions to diverge across runs. Finally, the importance of the interaction term `deg\_x\_bounds` (the product of the degeneracy index and the ratio of bounded variables) suggests that degeneracy is especially destabilizing when it occurs within a tightly constrained residual model. Overall, these results suggest that variability is a systematic property of the numerical and combinatorial structure exposed at the root node and point to directions for more informed solver tuning and fairer benchmarking.

\paragraph{Implications, Limitations, and Next Steps}\label{subsec:disc-implications}
The results have implications for how MIP benchmarks are used and how solvers are deployed. First, variability becomes a quantity that can be reported alongside standard aggregates: for each MIPLIB~2017 instance we obtain a scalar index $WVS$ (and, for SCIP, a time/gap/effort decomposition), which distinguishes genuinely stable instances from those whose outcomes depend strongly on the random seed. Reporting $WVS$ or at least flagging high-$WVS$ cases clarifies when an apparent performance gain is robust rather than sampling noise and lets researchers quantify how much variability a new branching rule, cut family, or heuristic introduces or removes. In practice, this does not mean running 100 seeds for every instance; rather, it suggests spending additional runs on the cases that appear highly variable, while treating low-variability instances with a lighter protocol. Second, the early-stage predictor suggests a simple triage mechanism. Because the model uses only static descriptors and root-level signals, it can be evaluated after presolve and the first LP relaxation: instances with high predicted $WVS$ are natural targets for larger budgets, extra seeds, or bet-and-run strategies, while low-$WVS$ instances appear safe for single-seed evaluation. The case of \textit{neos-5045105-creuse} illustrates this advantage: performance variability can provide a proof of optimality for a previously open instance if given multiple seeds. Third, the drivers highlighted by $WVS$ and the predictor point toward solver development. High predicted variability is tied to degenerate and dual-inactive root relaxations and the specific combinatorial structure remaining after presolve. This indicates that variability should be treated as an explicit tuning objective: parameter sets can be compared not only on speed and solved counts but also on their impact on $WVS$ for selected families of instances.

Because these conclusions rest on a finite sample, each $WVS$ value is only
an estimate of the underlying variability, and small differences in scores
or predictions should not be overinterpreted. Extending the study to additional solvers and versions, incorporating deeper tree-level signals, exploring task-specific weight choices, and coupling the forecaster to variability-aware policies are natural next steps.

\section{Conclusion}\label{sec:conclusion}

We introduced the Weighted Variability Score (WVS), a timeout-robust index that aggregates down-weighted dispersion in normalized runtime, final gap, LP iterations, and branch-and-bound nodes into a single, scale-free measure of solver instability. On a repeated-run study of MIPLIB~2017, WVS reveals a strongly right-skewed distribution with a substantial upper tail: many instances are essentially stable, but a nontrivial minority exhibit large run-to-run variability that is not a simple monotone function of instance scale or difficulty. Using only static descriptors and early-run features, a compact tree-based predictor with a gated tail uplift attains a mean absolute error of about $0.058$ and a Spearman correlation of about $0.65$ on a sealed test split, with better identification of high-variability instances than runtime-only proxies. Feature-importance analyses point consistently to root degeneracy, dual inactivity, and the numerical and combinatorial structure exposed at the root node as the main drivers. We hope that WVS and the accompanying predictor will serve as a practical toolkit for reporting, understanding, and ultimately controlling variability in mixed-integer programming.


%

%
%

\section*{Acknowledgments}
We thank Professor David Shmoys and the Center for Data Science for Enterprise and Society at Cornell University for supporting access to the cloud-computing resources used in this study.



\bibliographystyle{informs2014} 




\end{document}




\RUNAUTHOR{Adegbindin et al.}

\RUNTITLE{Performance Variability in MIP Solvers}

\TITLE{When Solvers Get Moody: Quantifying and Predicting Performance Variability in Mixed Integer Programming}

\ARTICLEAUTHORS{%
\AUTHOR{Mourchid Adegbindin}
\AFF{Department of Operations Research of Information Engineering,
Cornell University, \EMAIL{ma757@cornell.edu}}

\AUTHOR{Brenda Dietrich}
\AFF{Department of Operations Research of Information Engineering,
Cornell University, \EMAIL{bld34@cornell.edu}}

\AUTHOR{Jamol Pender}
\AFF{Department of Operations Research of Information Engineering,
Cornell University, \EMAIL{jjp274@cornell.edu}}

\AUTHOR{Andrea Lodi}
\AFF{Department of Operations Research of Information Engineering,
Cornell University, \EMAIL{al748@cornell.edu}}

} 







\ECSwitch 

\ECHead{Online Companion to:\\ When Solvers Get Moody: Quantifying and Predicting Performance Variability in Mixed Integer Programming}


%

\section{Selected Feature Definitions}

To aid interpretation of the feature-importance results, we define several
features that appear prominently in the analysis. Exact feature construction
is provided in the accompanying reproducibility code.

\begin{itemize}

    \item \texttt{frac\_zero\_rc}: Ratio of non-basic variables with zero
    reduced cost to all non-basic variables; used as a proxy for primal
    degeneracy.

    \item \texttt{N\_Bind\_Cons\_Zero\_Dual\_to\_N\_Cons\_ratio}: Number of
    binding constraints with zero dual value divided by the total number of
    constraints; used as a proxy for dual inactivity.

    \item \texttt{Presol\_Precedence}: Number of precedence-type constraints
    remaining after presolve.

    \item \texttt{Presol\_Nonzeroes}: Number of nonzeros in the constraint
    matrix after presolve.

    \item \texttt{LP\_Solving\_Time}: An elapsed-time proxy obtained from the
    first running-status statistics after presolving. At the instance level,
    this quantity is aggregated across randomized runs.

\end{itemize}

\section{Historical Comparison with MIPLIB 2010}
This section reports the 38 instances common to the MIPLIB~2010
variability study and our SCIP corpus. The legacy variability score (VS)
is the runtime coefficient of variation reported by
\citet{KochEtAl2011_MIPLIB2010}. The dagger marks the three legacy VS
values reported as affected by the 10-hour time limit used in that study.

\begin{longtable}{lrrr}
\caption{MIPLIB~2010 variability scores and corresponding variability measures for the 38 common instances. $^\dagger$ denotes a legacy VS value reported as affected by the 10-hour limit in the MIPLIB~2010 experiment. \label{tab:miplib2010-overlap}}\\
\hline
Instance & VS (2010) & Runtime CV & WVS \\
\hline
\endfirsthead
\multicolumn{4}{c}{\tablename\ \thetable\ - continued}\\
\hline
Instance & VS (2010) & Runtime CV & WVS \\
\hline
\endhead
\hline
\multicolumn{4}{r}{Continued on next page}\\
\endfoot
\hline
\endlastfoot
acc-tight5 & 0.72 & 0.609 & 0.479 \\
aflow40b & 0.39 & 0.255 & 0.191 \\
air04 & 0.14 & 0.182 & 0.192 \\
biella1 & 0.29 & 0.346 & 0.241 \\
bienst2 & 0.21 & 0.233 & 0.161 \\
binkar10\_1 & 0.50 & 0.208 & 0.131 \\
bley\_xl1 & 0.26 & 0.225 & 0.295 \\
core2536-691 & 0.67 & 0.806 & 0.437 \\
csched010 & 0.26 & 0.309 & 0.267 \\
danoint & 0.09 & 0.125 & 0.083 \\
eil33-2 & 0.15 & 0.100 & 0.088 \\
ex9 & 0.11 & 0.112 & 0.068 \\
lectsched-4-obj & 1.04 & 1.158 & 0.717 \\
map18 & 0.15 & 0.281 & 0.174 \\
mcsched & 0.16 & 0.181 & 0.138 \\
mine-166-5 & 0.24 & 0.189 & 0.146 \\
mine-90-10 & 0.75 & 0.487 & 0.309 \\
mzzv11 & 0.21 & 0.178 & 0.128 \\
neos-1396125 & 0.25 & 0.380 & 0.279 \\
neos-1601936 & 0.81$^\dagger$ & 0.355 & 0.328 \\
neos-686190 & 0.27 & 0.393 & 0.279 \\
neos18 & 0.39 & 0.238 & 0.171 \\
net12 & 0.33 & 0.644 & 0.385 \\
newdano & 0.21 & 0.216 & 0.229 \\
noswot & 0.34 & 0.368 & 0.249 \\
ns1208400 & 0.70 & 1.223 & 0.697 \\
ns1830653 & 0.42 & 0.361 & 0.231 \\
pg5\_34 & 0.13 & 0.247 & 0.155 \\
pw-myciel4 & 0.51 & 0.470 & 0.339 \\
qiu & 0.19 & 0.284 & 0.170 \\
rmatr100-p10 & 0.08 & 0.166 & 0.099 \\
rmatr100-p5 & 0.12 & 0.036 & 0.025 \\
rococoC10-001000 & 0.42 & 0.339 & 0.219 \\
roll3000 & 0.54$^\dagger$ & 0.284 & 0.200 \\
sp98ir & 0.20 & 0.309 & 0.203 \\
timtab1 & 0.43 & 0.218 & 0.143 \\
triptim1 & 0.48 & 0.690 & 0.472 \\
unitcal\_7 & 1.05$^\dagger$ & 0.194 & 0.163 \\
\end{longtable}

\section{Additional Reproducibility Details}

The computational workflow proceeds from raw solver outputs to a run-level
parsed data set, followed by instance-level aggregation, computation of the
component MCV values and WVS, and generation of the reported figures and
prediction results.

For the WVS calculation, runtimes are truncated to the interval
$[1,14400]$ seconds and divided by $14400$, solver-reported final relative
gaps are capped at one, and LP iterations and branch-and-bound nodes are used
as reported. Component standard deviations use denominator $N$. The fixed
WVS weights for time, gap, LP iterations, and nodes are respectively
$(0.65,0.05,0.20,0.10)$.

The complete run-level data and code used to reproduce the computational
results will be made publicly available upon publication.

%
%




\bibliographystyle{informs2014} 



